\documentclass[11pt,dvips]{article}
\usepackage{theorem}
\usepackage{groups}
\usepackage{epsfig}
\usepackage{amsfonts}
\usepackage{amssymb}
\usepackage{latexsym}
\usepackage{natbib}
\def\to{\longrightarrow}                
\def\span#1{\langle #1\rangle}          
\def\neutral{e}                         
\def\blank{\ }                          
\def\fullblank{\_}                      
\def\field#1{GF(#1)}                    
\def\det#1{\mathrm{det}\ #1}            
\def\dpr#1#2{#1\cdot #2}                
\def\vprsign{\times}                    
\def\vpr#1#2{#1\vprsign #2}             
\def\aut#1{\mathrm{Aut}(#1)}            
\def\centre#1{Z(#1)}                    
\def\order#1{|#1|}                      
\def\cyclic#1{C_{#1}}                   
\def\vm#1#2#3#4{\left(\begin{array}{cc}#1 & #2\\ #3 & #4 \end{array}\right)}
\def\onezorn#1{M(#1)}                   
\def\paige#1{M^*(#1)}                   
\def\inv#1#2{[#1,\,#2]}                 
\def\iinv#1#2#3{\inv{#1}{#2}_{#3}}      
\def\tri#1#2#3{\{#1,\,#2\}_{#3}}        
\def\weight#1{\mathrm{w}(#1)}           
\def\commutator#1#2{[#1,\,#2]}          
\def\normletter{N}                      
\def\norm#1{\normletter(#1)}            
\def\bilin#1#2{\normletter(#1,\,#2)}    
\def\octoletter{\mathbb O}              
\def\octo#1{\octoletter(#1)}            
\def\diagm#1{\widehat{#1}}              
\def\dswitch{\partial}                  
\def\lie#1{\mathrm{Lie}(#1)}            

\begin{document}

\maketitle{Chevalley Groups of Type $G_2$ as Automorphism Groups of Loops}{Petr
Vojt\v echovsk\'y}{Department of Mathematics, Iowa State University, Ames,
Iowa, 50011, U.S.A.}

\begin{abstract}
Let $\paige{q}$ be the unique nonassociative finite simple Moufang loop
constructed over $\field{q}$. We prove that $\aut{\paige{2}}$ is the Chevalley
group $G_2(2)$, by extending multiplicative automorphism of $\paige{2}$ into
linear automorphisms of the unique split octonion algebra over $\field{2}$.
Many of our auxiliary results apply in the general case. In the course of the
proof we show that every element of a split octonion algebra can be written as
a sum of two elements of norm one.
\end{abstract}

\section{Composition Algebras and Paige Loops}\label{Sc:CAPL}

\noindent Let $C$ be a finite-dimensional vector space over a field $k$,
equipped with a quadratic form $\normletter:C\to k$ and a multiplicative
operation $\cdot$. Following \cite{SpringerVeldkamp}, we say that
$C=(C,\,\normletter,\,+,\,\cdot)$ is a \emph{composition algebra} if
$(C,\,+,\cdot)$ is a \emph{nonassociative} ring with identity element
$\neutral$, $\normletter$ is nondegenerate, and $\norm{u\cdot v} =
\norm{u}\norm{v}$ is satisfied for every $u$, $v\in C$. The bilinear form
associated with $\normletter$ will also be denoted by $\normletter$. Recall
/that $\normletter: C\times C\to k$ is defined by
$\bilin{u}{v}=\norm{u+v}-\norm{u}-\norm{v}$. Write $u\perp v$ if $\bilin{u}{v}
= 0$, and set $u^\perp =\{v\in C;\; u\perp v\}$.

The standard $8$-dimensional real Cayley algebra $\octoletter$ constructed by
the \emph{Cayley-Dickson process} (or \emph{doubling} \cite{SpringerVeldkamp})
is the best known nonassociative composition algebra. There is a remarkably
compact way of constructing $\octoletter$ that avoids the iterative
Cayley-Dickson process. As in \cite{Coxeter}, let $B=\{\neutral=e_0$, $e_1$,
$\dots$, $e_7\}$ be a basis whose vectors are multiplied according to
\begin{equation}\label{Eq:Compact}
\begin{array}{l}
    e_r^2=-1,\;\;e_{r+7}=e_r,\;\;e_re_s=-e_se_r,\;\;
    e_{r+1}e_{r+3}=e_{r+2}e_{r+6}=e_{r+4}e_{r+5}=e_r,
\end{array}
\end{equation}
for $r$, $s\in\{1,\,\dots,\,7\}$, $r\ne s$. (Alternatively, see \cite[p.\ 122]
{ConwaySloane}.) The norm $\norm{u}$ of a vector
$u=\sum_{i=0}^7a_ie_i\in\octoletter$ is given by $\sum_{i=0}^7a_i^2$.

Importantly, all the \emph{structural constants} $\gamma_{ijk}$, defined by
$e_i\cdot e_j=\sum_{k=0}^7 \gamma_{ijk}e_k$, are equal to $\pm 1$, and
therefore the construction can be imitated over any field $k$. For
$k=\field{q}$ of odd characteristic, let us denote the ensuing algebra by
$\octo{q}$. When $q$ is even, the above construction does not yield a
composition algebra.

The following facts about composition algebras can be found in
\cite{SpringerVeldkamp}. Every non-trivial composition algebra $C$ has
dimension $2$, $4$ or $8$, and we speak of a \emph{complex}, \emph{quaternion}
or \emph{octonion} algebra, respectively. We say that $C$ is a \emph{division
algebra} if it has no zero divisors, else $C$ is called \emph{split}. There can
be many non-isomorphic octonion algebras over a given field. Exactly one of
them is guaranteed to be split. Moreover, when $k$ is finite, all octonion
algebras over $k$ are isomorphic (and thus split). Let $\octo{q}$ be the unique
octonion algebra constructed over $\field{q}$.

All composition algebras satisfy the so-called \emph{Moufang identities}
\begin{equation}\label{Eq:MoufangIdentities}
    (xy)(zx)=x((yz)x),\;\;
    x(y(xz))=((xy)x)z,\;\;
    x(y(zy))=((xy)z)y.
\end{equation}
These identities are the essence of Moufang loops, undoubtedly the most
investigated variety of nonassociative loops. More precisely, a quasigroup
$(L,\,\cdot)$ is a \emph{Moufang loop} if it possesses a neutral element
$\neutral$ and satisfies one (and hence all) of the Moufang identities
$(\ref{Eq:MoufangIdentities})$. We refer the reader to \cite{Pflugfelder} for
the basic properties of loops and Moufang loops in particular. Briefly, every
element $x$ of a Moufang loop $L$ has a both-sided inverse $x^{-1}$, and a
subloop $\span{x,\,y,\,z}$ of $L$ generated by $x$, $y$ and $z$ is a group if
and only if $x$, $y$ and $z$ associate. Specifically, every two-generated
subloop of $L$ is a group.

Paige \cite{Paige} constructed one nonassociative finite simple Moufang loop
for every finite field $\field{q}$. Liebeck \cite{Liebeck} used the
classification of finite simple groups in order to prove that there are no
other nonassociative finite simple Moufang loops. Reflecting the current trend
in loop theory, we will call these loops \emph{Paige loops}, and we denote the
unique Paige loop constructed over $\field{q}$ by $\paige{q}$.

The relation between $\octo{q}$ and $\paige{q}$ is as follows. Let
$\onezorn{q}$ be the set of all elements of $\octo{q}$ of norm one. Then
$\onezorn{q}$ is a Moufang loop with center
$\centre{\onezorn{q}}=\{\neutral,\,-\neutral\}$, and
$\onezorn{q}/\centre{\onezorn{q}}$ is isomorphic to $\paige{q}$. Note that
$\onezorn{q}=\paige{q}$ in characteristic $2$.

Historically, all split octonion algebras and Paige loops were constructed by
Zorn \cite{Zorn} and Paige without reference to doubling.  Given a field $k$,
consider the \emph{vector matrix algebra} consisting of all \emph{vector
matrices}
\begin{displaymath}
    x=\vm{a}{\alpha}{\beta}{b},
\end{displaymath}
where $a$, $b\in k$, $\alpha$, $\beta\in k^3$, addition is defined
entry-wise, and multiplication by
\begin{displaymath}
    \vm{a}{\alpha}{\beta}{b}\vm{c}{\gamma}{\delta}{d}
    =\vm{ac+\dpr{\alpha}{\delta}}{a\gamma+d\alpha-\vpr{\beta}{\delta}}
        {c\beta+b\delta+\vpr{\alpha}{\gamma}}{\dpr{\beta}{\gamma}+bd}.
\end{displaymath}
Here, $\dpr{\alpha}{\beta}$ (resp.\ $\vpr{\alpha}{\beta}$) is the standard dot
product (resp.\ vector product) of $\alpha$ and $\beta$. Use
$\det{x}=ab-\dpr{\alpha}{\beta}$ as a norm to obtain an octonion algebra. In
fact, this is exactly the unique split octonion algebra over $k$. The identity
element is
\begin{displaymath}
    \neutral=\vm{1}{(0,0,0)}{(0,0,0)}{1},
\end{displaymath}
and, when $\norm{x}=\det{x}\ne 0$, we have
\begin{displaymath}
    x^{-1}=\vm{b}{-\alpha}{-\beta}{a}.
\end{displaymath}

The purpose of this paper is to initiate the investigation of automorphism
groups of Paige loops. We will extend automorphisms of $\paige{2}$ into
automorphisms of $\octo{2}$ to prove that $\aut{\paige{2}}$ is the exceptional
Chevalley group $G_2(2)$. We also present several results for the general case.
Since $G_2(q)\le\aut{\onezorn{q}}$, it is reasonable to expect that equality
holds whenever $q$ is prime. See Acknowledgement for more details. It is fun to
watch how much information about the (boring) additive structure of a
composition algebra can be obtained from the multiplication alone (cf. Lemma
\ref{Lm:MultToAdd}).

\section{Multiplication versus Addition}

\noindent Perhaps the single most important feature of composition algebras is
the existence of the minimal equation (cf. \cite[Prop
1.2.3]{SpringerVeldkamp}). Namely, every element $x\in C$ satisfies
\begin{equation}\label{Eq:Minimal}
    x^2 - \bilin{x}{\neutral}x + \norm{x}\neutral=0.
\end{equation}
Furthermore, $(\ref{Eq:Minimal})$ is the minimal equation for $x$ when $x$ is
not a scalar multiple of $\neutral$.

\begin{lemma}\label{Lm:NormZero}
Let $C$ be a composition algebra, $x$, $y\in C$. Then
\begin{equation}\label{Eq:SpecialBilin}
    \bilin{xy}{y}=\bilin{x}{\neutral}\norm{y}.
\end{equation}
When $\norm{y}\ne 0$, we have
\begin{equation}\label{Eq:MinimalXY}
    (xy^{-1})^2-\bilin{x}{y}\norm{y}^{-1}xy^{-1}+\norm{xy^{-1}}\neutral=0.
\end{equation}
In particular,
\begin{equation}
    (xy^{-1})^2-\bilin{x}{y}xy^{-1}+\neutral=0 \label{Eq:SpecialMinXY}
\end{equation}
whenever $\norm{x}=\norm{y}=1$. In such a case, $(xy^{-1})^2=-\neutral$ if and
only if $x\perp y$.
\end{lemma}
\begin{proof}
We have $\bilin{xy}{y}=\norm{xy+y}-\norm{x}\norm{y}-\norm{y}$, and
$\norm{xy+y}=\norm{x+\neutral}\norm{y}=(\bilin{x}{\neutral}+
\norm{x}+\norm{\neutral})\norm{y}=
\bilin{x}{\neutral}\norm{y}+\norm{x}\norm{y}+\norm{y}$. Equation
$(\ref{Eq:SpecialBilin})$ follows.

Substitute $xy^{-1}$ for $x$ into $(\ref{Eq:SpecialBilin})$ to obtain
$\bilin{x}{y}=\bilin{xy^{-1}}{\neutral}\norm{y}$. The minimal equation
\begin{displaymath}
    (xy^{-1})^2-\bilin{xy^{-1}}{\neutral}xy^{-1}+\norm{xy^{-1}}\neutral=0
\end{displaymath}
for $xy^{-1}$ can then be written as $(\ref{Eq:MinimalXY})$, provided
$\norm{y}\ne 0$. The rest is easy.
\end{proof}

We leave the proof of the following Lemma to the reader.

\begin{lemma}\label{Lm:TwoAndThree}
Let
\begin{displaymath}
    x=\vm{a}{\alpha}{\beta}{b}
\end{displaymath}
be an element of a composition algebra $C$ satisfying $\norm{x}=1$. Then
\begin{enumerate}
\item[\emph{(i)}] $x^2=\neutral$ if and only if $(\alpha,\,\beta)=(0,\,0)$
    and $a=b\in\{1,\,-1\}$,
\item[\emph{(ii)}] $x^2=-\neutral$ if and only if $((\alpha,\,\beta)=
    (0,\,0)$, $b=a^{-1}$, and $a^2=-1)$ or $((\alpha,\,\beta)\ne(0,\,0)$ and
    $b=-a)$,
\item[\emph{(iii)}] $x^3=\neutral$ if and only if $((\alpha,\,\beta)=
    (0,\,0)$, $b=a^{-1}$, and $a^3=1)$ or $((\alpha,\,\beta)\ne(0,\,0)$ and
    $b=-1-a)$,
\item[\emph{(iv)}] $x^3=-\neutral$ if and only if $((\alpha,\,\beta)=
    (0,\,0)$, $b=a^{-1}$, and $a^3=-1)$ or $((\alpha,\,\beta)\ne(0,\,0)$ and
    $b=1-a)$.
\end{enumerate}
\end{lemma}

Let us denote the multiplicative order of $x$ by $\order{x}$.

\begin{lemma}\label{Lm:NormOne}
Let $C$ be a composition algebra. Assume that $x$, $y\in C$ satisfy
$\norm{x}=\norm{y}=1$, $x\ne y$. The following conditions are equivalent:
\begin{enumerate}
\item[\emph{(i)}] $\order{xy^{-1}}=3$,
\item[\emph{(ii)}] $(xy^{-1})^2+xy^{-1}+\neutral=0$,
\item[\emph{(iii)}] $\bilin{x}{y}=-1$,
\item[\emph{(iv)}] $\norm{x+y}=1$.
\end{enumerate}
\end{lemma}
\begin{proof}
The equivalence of (ii) and (iii) follows from the uniqueness of the minimal
equation $(\ref{Eq:Minimal})$. Condition (iii) is equivalent to (iv) since
$\norm{x}=\norm{y}=1$. It suffices to prove the equivalence of (i) and (ii).

As $(a^3-\neutral)=(a-\neutral)(a^2+a+\neutral)$, there is nothing to prove
when $C$ has no zero divisors. The implication (ii) $\Rightarrow$ (i) is
obviously true in any (composition) algebra. Let us prove (i) $\Rightarrow$
(ii) for a split octonion algebra $C$. Assume that $\order{xy^{-1}}=3$, $x\ne
y$, and
\begin{displaymath}
    x=\vm{a}{\alpha}{\beta}{b},\;\; y=\vm{c}{\gamma}{\delta}{d}.
\end{displaymath}
We prove that $\norm{x+y}=1$. Direct computation yields $\norm{x+y}=2+r+s$,
where $r=ad-\dpr{\alpha}{\delta}$, $s=bc-\dpr{\beta}{\gamma}$. Also,
\begin{displaymath}
    xy^{-1}=\vm{r}{\varepsilon}{\varphi}{s}
\end{displaymath}
for some $\varepsilon$, $\varphi\in k^3$. Since $(xy^{-1})^3=\neutral$, we have
either $((\varepsilon,\,\varphi)=(0,\,0)$, $s=r^{-1}$, and $r^3=1)$, or
$((\varepsilon,\,\varphi)\ne(0,\,0)$, and $s=-1-r)$, by Lemma
\ref{Lm:TwoAndThree}. If the latter is true, we immediately get $\norm{x+y}=1$.
Assume the former is true. Then $r+s=r+r^{-1}$. Also, $r^3=1$ implies $r=1$ or
$r^2+r+1=0$. But $r=1$ leads to $x=y$, a contradiction. Therefore $r^2+r+1=0$,
i.e., $r+r^{-1}=-1$, and we get $\norm{x+y}=1$ again.
\end{proof}

There is a strong relation between the additive and multiplicative structures
in composition algebras.

\begin{lemma}\label{Lm:MultToAdd}
Let $C$, $x$, and $y$ be as in Lemma $\ref{Lm:NormOne}$. Then $\norm{x+y}=1$ if
and only if $x+y=-xy^{-1}x$.
\end{lemma}
\begin{proof}
The indirect implication is trivial. Assume that $\norm{x+y}=1$. Then
$(xy^{-1})^2 + xy^{-1}+\neutral=0$, and $(xy^{-1})^3=\neutral$. Thus
$yx^{-1}=(xy^{-1})^2=-xy^{-1}-\neutral$. Multiplying this equality on the right
by $x$ yields $y=-xy^{-1}x-x$.
\end{proof}

\section{Doubling Triples}

\noindent Any composition algebra $C$ can be constructed from the underlying
field $k$ in three steps. Proposition $1.5.1$ and Lemma $1.6.1$ of
\cite{SpringerVeldkamp} tell us how to do it. Imitating these results, we say
that a triple $(a,\, b,\, c)\in C^3$ is a \emph{doubling triple} if
$\norm{a}\ne 0$, $\norm{b}\ne 0$, $\norm{c}\ne 0$, $b\in\neutral^\perp\cap
a^\perp$, $c\in\neutral^\perp\cap a^\perp\cap b^\perp\cap (ab)^\perp$,
$a\in\neutral^\perp$ (resp.\ $a\not\in\neutral^\perp$) and the characteristic
of $k$ is odd (resp.\ even). Then $B=\{\neutral$, $a$, $b$, $ab$, $c$, $ac$,
$bc$, $(ab)c\}$ is a basis for $C$.

Construction $(\ref{Eq:Compact})$ immediately shows that there is a doubling
triple with $\norm{a}=\norm{b}=\norm{c}=1$ when $k$ is of odd characteristic.
Such a doubling triple exists for every split octonion algebra in even
characteristic, too.

\begin{lemma}\label{Lm:DoublingTriple}
Let $k$ be a field of even characteristic. Then
\begin{displaymath}
     \left( \vm{0}{(1,0,0)}{(1,0,0)}{1},\;\;
            \vm{0}{(0,1,0)}{(0,1,0)}{0},\;\;
            \vm{0}{(0,0,1)}{(0,0,1)}{0}
    \right)
\end{displaymath}
is a doubling triple consisting of elements of norm one.
\end{lemma}
\begin{proof}
Straightforward computation.
\end{proof}

Doubling triples can be used to induce automorphisms. By an automorphism of a
composition algebra $C$ we mean a \emph{linear} automorphism, i.e., a bijection
$f:C\to C$ satisfying $f(u+v)=f(u)+f(v)$, $f(\lambda u)=\lambda f(u)$, and
$f(u\cdot v)=f(u)\cdot f(v)$ for every $u$, $v\in C$, $\lambda\in k$. By
\cite[Thm 1.7.1 and Cor 1.7.2]{SpringerVeldkamp}, every such automorphism is an
isometry (i.e., $\norm{f(u)}=\norm{u}$), and vice versa. Springer and Veldkamp
\cite[Ch. 2]{SpringerVeldkamp} use algebraic groups to show that
$\aut{\octo{q}}$ is the exceptional group $G_2(q)$, and more.

\begin{prop}\label{Pr:ExtendInSteps}
Let $(a,\,b,\,c)$, $(a',\,b',\,c')$ be two doubling triples of a composition
algebra $C$. Then there is an automorphism of $C$ mapping $(a,\,b,\,c)$ onto
$(a',\,b',\,c')$ if and only if $\norm{x}=\norm{x'}$, for $x=a$, $b$, $c$.
\end{prop}
\begin{proof} Let $k$ be the underlying field.
The necessity is obvious since every automorphism is an isometry. Now for the
sufficiency. Let $A=k\neutral\oplus ka$, $B=A\oplus Ab$, $C=B\oplus Bc$, and
similarly for $A'$, $B'$, $C'=C$. Define $\psi_X:X\to X'=\psi_X(X)$ (for $X=A$,
$B$, $C$) by
\begin{eqnarray*}
    \psi_A(x+ya)&=&x+ya' \;\;\;(x,\,y\in k),\\
    \psi_B(x+yb)&=&\psi_A(x)+\psi_A(y)b' \;\;\; (x,\,y\in A),\\
    \psi_C(x+yc)&=&\psi_B(x)+\psi_B(y)c' \;\;\;(x,\,y\in B).
\end{eqnarray*}
All maps $\psi_X$ are clearly linear, and it is not hard to see that they are
also multiplicative. (One has to use the assumption $\norm{x}=\norm{x'}$.)
Since $a'$, $b'$, $c'$ generate $C'=C$, $\psi_C$ is the automorphism we are
looking for.
\end{proof}

\section{Restrictions and Extensions of Automorphisms}\label{Sc:RE}

\noindent The restriction of $h\in\aut{\octo{q}}$ onto the loop $\onezorn{q}$
is a (multiplicative) automorphism. Moreover, two distinct automorphisms of
$\octo{q}$ differ on $\onezorn{q}$, because there is a basis for $C$ consisting
of unit vectors (cf.\ construction (\ref{Eq:Compact}) and Lemma
\ref{Lm:DoublingTriple}).

We would like to emphasize at this point how far are the metric properties of
$\normletter$ from our intuitive understanding of (real) norms. Theorem
\ref{Th:SumsOfTwo} is not required for the rest of the paper, but is certainly
of interest in its own right.

\begin{thm}\label{Th:SumsOfTwo}
Every element of a split octonion algebra $C$ is a sum of two elements of norm
one.
\end{thm}
\begin{proof}
We identify $C$ with the vector matrix algebra over $k$, where the norm is
given by the determinant. Let
\begin{displaymath}
    x=\vm{a}{\alpha}{\beta}{b}
\end{displaymath}
be an element of $C$. First assume that $\beta\ne 0$. Note that for every
$\lambda\in k$ there is $\gamma\in k^3$ such that
$\dpr{\gamma}{\beta}=\lambda$. Pick $\gamma\in k^3$ so that
$\dpr{\gamma}{\beta}=a+b-ab+\dpr{\alpha}{\beta}$. Then choose
$\delta\in\gamma^\perp\cap\alpha^\perp\ne\emptyset$. (As usual,
$\alpha\perp\delta$ if and only if $\dpr{\alpha}{\delta}=0$.) This choice
guarantees that $(a-1)(b-1)-\dpr{(\alpha-\gamma)}{(\beta-\delta)}
=ab-a-b+1-\dpr{\alpha}{\beta}+\dpr{\gamma}{\beta}=1$. Thus
\begin{displaymath}
    \vm{a}{\alpha}{\beta}{b}=\vm{1}{\gamma}{\delta}{1}
        +\vm{a-1}{\alpha-\gamma}{\beta-\delta}{b-1}
\end{displaymath}
is the desired decomposition of $x$ into a sum of two elements of norm $1$.
Note that the above procedure works even for $\alpha=0$.

Now assume that $\beta=0$. If $\alpha\ne 0$, we use a symmetrical argument as
before to decompose $x$. It remains to discuss the case when $\alpha=\beta=0$.
Then the equality
\begin{displaymath}
    \vm{a}{0}{0}{b}=\vm{a}{(1,0,0)}{(-1,0,0)}{0}+\vm{0}{(-1,0,0)}{(1,0,0)}{b}
\end{displaymath}
does the job.
\end{proof}

We now know that $G_2(q)$ is a subgroup of $\aut{\onezorn{q}}$. Let us consider
the extension problem. Pick an automorphism $g$ of the (not necessarily simple)
Moufang loop $\onezorn{q}$. The ultimate goal is to construct
$h\in\aut{\octo{q}}$ such that $h\restriction \onezorn{q}=g$. If this can be
done, we immediately conclude that $\aut{\onezorn{q}}=G_2(q)$ for every $q$. We
like to think of the problem as a notion ``orthogonal'' to Witt's lemma.
Roughly speaking, Witt's lemma deals with extensions of partial isometries from
subspaces onto finite-dimensional vector spaces, whereas we are attempting to
extend a multiplicative, norm-preserving map from the first shell $\onezorn{q}$
into an automorphism ($=$ isometry) of $\octo{q}$. Naturally, $g$ is not linear
because $\onezorn{q}$ is not even closed under addition. However, the analogy
with Witt's lemma will become more apparent once we prove that $g$ \emph{is},
in a sense, additive (cf. Proposition \ref{Pr:Additivity}).

\begin{lemma}\label{Lm:DoublingTriples}
Let $g\in\aut{\onezorn{q}}$, and let $(a,\,b,\,c)$ be a doubling triple for
$\octo{q}$ with $\norm{a}=\norm{b}=\norm{c}=1$. Then $(g(a),\,g(b),\,g(c))$ is
a doubling triple $($with $\norm{g(a)}=\norm{g(b)}=\norm{g(c)}=1)$.
\end{lemma}
\begin{proof}
Since $(a,\,b,\,c)$ is a doubling triple, we have $b\in\neutral^\perp\cap
a^\perp$, $c\in\neutral^\perp\cap a^\perp\cap b^\perp$. Moreover,
$a\in\neutral^\perp$ (resp. $a\not\in\neutral^\perp$) if $q$ is odd (resp.\
even). By Lemmas \ref{Lm:NormZero} and \ref{Lm:NormOne}, this is equivalent to
$b^2=c^2=(ab^{-1})^2=(ac^{-1})^2=(bc^{-1})^2 = ((ab)c^{-1})^2 = -\neutral$, and
$a^2=-\neutral$ (resp.\ $|a|=3$). Because $g\in\aut{\onezorn{q}}$, we have
$g(b)^2=g(c)^2 = (g(a)g(b)^{-1})^2=(g(a)g(c)^{-1})^2=(g(b)g(c)^{-1})^2 =
((g(a)g(b))g(c)^{-1})^2 = -\neutral$ and $g(a)^2=-\neutral$ (resp.\
$|g(a)|=3$). Another application of Lemmas \ref{Lm:NormZero} and
\ref{Lm:NormOne} shows that $(g(a)$, $g(b)$, $g(c)$) is a doubling triple.
\end{proof}

In particular, the mapping $h=\psi_{\octo{q}}$ constructed from $g$ and
$(a,\,b,\,c)$ by Proposition \ref{Pr:ExtendInSteps} is an automorphism of
$\octo{q}$ satisfying $\psi(x)=g(x)$, for $x=a$, $b$, $c$.

\begin{re}
This extension $h$ can be obtained in another way when the characteristic is
odd. Namely, construct $\octo{q}$ as in section \ref{Sc:CAPL}, and define
$h:\octo{q}\to\octo{q}$ by
\begin{displaymath}
    h(\sum_{i=0}^7 a_ie_i) = \sum_{i=0}^7 a_ig(e_i).
\end{displaymath}
Obviously, $h$ is linear. For fixed $i$, $j$, only one of the $8$ structural
constants $\gamma_{ijk}$ is nonzero, and it is equal to $\pm 1$. Using
linearity of $h$, it is therefore easy to check that $h$ is multiplicative.
\end{re}

By the construction, $h$ coincides with $g$ on a basis $B$. However, we do not
know whether $h$ is an extension of $g$. The fact that $h\restriction B =
g\restriction B$ does not guarantee that $h\restriction \onezorn{q} = g$, since
$B$ does not need to generate $\onezorn{q}$ by multiplication. Interestingly
enough, it seems to never be the case! The key to answering these questions is
to look at the additive properties of $g$.

\section{Automorphisms of Finite Octonion Algebras}\label{Sc:Autos}

\noindent We have entered a more technical part of the paper. In this section,
we construct a family of automorphisms of $\octo{q}$.

Let $k=\field{q}$, and let $\lie{q}$ be the three-dimensional Lie algebra $k^3$
with vector product $\times$ playing the role of a Lie bracket. A linear
transformation $f:\lie{q}\to\lie{q}$ belongs to $\aut{\lie{q}}$ if and only if
$f(\vpr{\alpha}{\beta})=\vpr{f(\alpha)}{f(\beta)}$ is satisfied for every
$\alpha$, $\beta\in k^3$. We say that a linear transformation $f$ is
\emph{orthogonal} if $\dpr{f(\alpha)}{f(\beta)}=\dpr{\alpha}{\beta}$ for every
$\alpha$, $\beta\in k^3$.

\begin{lemma}\label{Lm:SomeAutomorphisms}
For a non-singular orthogonal linear transformation $f:k^3\to k^3$, let
$\diagm{f}:\octo{q}\to\octo{q}$ be the mapping
\begin{displaymath}
    \diagm{f}\vm{a}{\alpha}{\beta}{b}=\vm{a}{f(\alpha)}{f(\beta)}{b}.
\end{displaymath}
Then $\diagm{f}\in\aut{\octo{q}}$ if and only if $f\in\aut{\lie{q}}$.
\end{lemma}
\begin{proof}
The map $\diagm{f}$ is clearly linear and preserves the norm. Since $f$ is
one-to-one, so is $\diagm{f}$. We have
\begin{eqnarray*}
    &&\diagm{f}\vm{a}{\alpha}{\beta}{b}\diagm{f}\vm{c}{\gamma}{\delta}{d}\\
    &&=\vm{ac+\dpr{f(\alpha)}{f(\delta)}}
        {af(\gamma)+df(\alpha)-\vpr{f(\beta)}{f(\delta)}}
        {cf(\beta)+bf(\delta)+\vpr{f(\alpha)}{f(\gamma)}}
        {\dpr{f(\beta)}{f(\gamma)}+bd}.
\end{eqnarray*}
On the other hand,
\begin{displaymath}
    \diagm{f}\left(\vm{a}{\alpha}{\beta}{b}
        \vm{c}{\gamma}{\delta}{d}\right)
    =\vm{ac+\dpr{\alpha}{\delta}}
         {f(a\gamma+d\alpha-\vpr{\beta}{\delta})}
         {f(c\beta+b\delta+\vpr{\alpha}{\gamma})}
         {\dpr{\beta}{\gamma}+bd}.
\end{displaymath}
Sufficiency is now obvious, and necessity follows by specializing the
elements $a$, $b$, $c$, $d$, $\alpha$, $\beta$, $\gamma$, $\delta$.
\end{proof}

For a map $f$, let $-f$ be the map \emph{opposite} to $f$, i.e.,
$(-f)(u)=-(f(u))$. Also, for a permutation $\pi\in S_3$, consider $\pi$ as a
linear transformation on $k^3$ defined by
\begin{displaymath}
    \pi(\alpha_1,\,\alpha_2,\,\alpha_3)=
        (\alpha_{\pi(1)},\,\alpha_{\pi(2)},\,\alpha_{\pi(3)}).
\end{displaymath}
Apparently, $-S_3=\{-\pi;\;\pi\in S_3\}$ is a set of non-singular orthogonal
linear transformations.

\begin{lemma}\label{Lm:ShuffleCoords}
$\diagm{-\pi}\in\aut{\octo{q}}$ for every $\pi\in S_3$.
\end{lemma}
\begin{proof}
Let $\pi\in S_3$ be the transposition interchanging $1$ and $2$, and let
$\alpha$, $\beta\in k^3$. Then
\begin{displaymath}
    \pi(\alpha\times\beta)=(\alpha_3\beta_1-\alpha_1\beta_3,\,
                            \alpha_2\beta_3-\alpha_3\beta_2,\,
                            \alpha_1\beta_2-\alpha_2\beta_1),
\end{displaymath}
and
\begin{displaymath}
    \pi(\alpha)\times\pi(\beta)=(\alpha_1\beta_3-\alpha_3\beta_1,\,
                                 \alpha_3\beta_2-\alpha_2\beta_3,\,
                                 \alpha_2\beta_1-\alpha_1\beta_2).
\end{displaymath}
Hence $-\pi(\alpha\times\beta)=
\pi(\alpha)\times\pi(\beta)=(-\pi)(\alpha)\times(-\pi)(\beta)$. Thanks to the
symmetry of $S_3$, we have shown that $-\pi\in\aut{\lie{q}}$ for every $\pi\in
S_3$. The rest follows from Lemma \ref{Lm:SomeAutomorphisms}.
\end{proof}

Observe there is another automorphism when $q$ is even:

\begin{lemma}\label{Lm:DiagonalSwitch}
Define $\dswitch:\octo{q}\to\octo{q}$ by
\begin{displaymath}
    \dswitch\vm{a}{\alpha}{\beta}{b}=\vm{b}{\beta}{\alpha}{a}.
\end{displaymath}
Then $\dswitch\in\aut{\octo{q}}$ if and only if $q=2^n$.
\end{lemma}

Finally, we look at conjugations. Let $L$ be a Moufang loop. For $x\in L$,
define the conjugation $T_x:L\to L$ by $T_x(y)=x^{-1}yx$, where $x^{-1}yx$ is
unambiguous thanks to the properties of $L$. Not every conjugation of $L$ is an
automorphism. By \cite[Thm IV.1.6]{Pflugfelder}, $T_x\in\aut{L}$ if
$x^3=\neutral$. (And it is not difficult to show that $x^3=\neutral$ is also a
necessary condition, provided $L$ is simple.)

\section{Transitivity of the Natural Action}\label{Sc:Trans}

\noindent We take advantage of the automorphisms defined in section
\ref{Sc:Autos}, and investigate the natural action of $\aut{\paige{2}}$ on
$\paige{2}=\onezorn{2}$. The lattice of subloops of $\paige{2}$ was fully
described in \cite{PetrPhD}. Here, we only focus on the action of
$\aut{\paige{2}}$ on involutions and on subgroups of $\paige{2}$ isomorphic to
$V_4$.

Once again, identify $\octo{2}$ with the vector matrix algebra.

By Lemma \ref{Lm:TwoAndThree}, $\paige{2}$ contains only elements of order $1$,
$2$ and $3$. Furthermore, every involution $x\in\paige{2}$ is of the form
\begin{displaymath}
    \vm{a}{\alpha}{\beta}{a},
\end{displaymath}
for some $a\in\{0,\,1\}$ and $\alpha$, $\beta\in k^3$. In order to linearize
our notation, we write $x=\inv{\alpha}{\beta}$ when the value of $a$ is clear
from $\alpha$, $\beta$, or when it is not important; and
$x=\iinv{\alpha}{\beta}{a}$ otherwise.

Similarly, every element $x\in\paige{2}$ of order $3$ is of the form
\begin{displaymath}
    \vm{a}{\alpha}{\beta}{1+a}.
\end{displaymath}
To save space, we write $x=\tri{\alpha}{\beta}{a}$.

We will sometimes leave out commas and parentheses. Thus, both $110$ and
$(110)$ stand for $(1,1,0)$. The commutator of $x$ and $y$ will be denoted by
$\commutator{x}{y}$.

\begin{prop}\label{Pr:Involutions}
Let $x=\iinv{\alpha}{\beta}{a}$, $y=\iinv{\gamma}{\delta}{b}$ be two
involutions of $\paige{2}$, $x\ne y$. Then:
\begin{enumerate}
\item[\emph{(i)}]
    $\commutator{x}{y}=\neutral$ if and only if $\order{xy}=2$ if and only if
    $\span{x,\,y}\cong V_4$ if and only if
    $\dpr{\alpha}{\delta}=\dpr{\beta}{\gamma}$.
\item[\emph{(ii)}]
    $\commutator{x}{y}\ne\neutral$ if and only if $\order{xy}=3$ if and only if
    $\span{x,\,y}\cong S_3$ if and only if
    $\dpr{\alpha}{\delta}\ne\dpr{\beta}{\gamma}$.
\item[\emph{(iii)}]
    $x$ is contained in a subgroup isomorphic to $S_3$,
\item[\emph{(iv)}]
    every subgroup of $\paige{2}$ isomorphic to $S_3$ contains an involution of the form
    $\iinv{\blank}{\blank}{0}$.
\end{enumerate}
\end{prop}
\begin{proof}
The involution $x$ commutes with $y$ if and only if $\order{xy}=2$. Since
\begin{displaymath}
    xy=\vm{ab+\dpr{\alpha}{\delta}}{\fullblank}{\fullblank}
        {ab+\dpr{\beta}{\gamma}},
\end{displaymath}
parts (i) and (ii) follow.

Given $x=\iinv{\alpha}{\beta}{a}$, pick $\delta\in\alpha^\perp$,
$\gamma\not\in\beta^\perp$, and choose $b\in\{0,\,1\}$ so that
$y=\iinv{\gamma}{\delta}{b}\in\paige{2}$. Then $\span{x,\,y}\cong S_3$, and
(iii) is proved.

Let $G\le\paige{2}$, $G\cong S_3$. Let $x=\iinv{\alpha}{\beta}{1}$,
$y=\iinv{\gamma}{\delta}{1}\in G$, $x\ne y$. Then
\begin{displaymath}
    xy=\vm{1+\dpr{\alpha}{\delta}}{\alpha+\gamma+\vpr{\beta}{\delta}}
        {\beta+\delta+\vpr{\alpha}{\gamma}}{1+\dpr{\beta}{\gamma}}.
\end{displaymath}
Since $\order{xy}=3$, we have $\dpr{\alpha}{\delta}\ne\dpr{\beta}{\gamma}$. In
other words, $\dpr{\alpha}{\delta}+\dpr{\beta}{\gamma}=1$. Then the third
involution $xyx\in G$ equals
\begin{displaymath}
    \vm{1+\dpr{\alpha}{\delta}+\dpr{(\alpha+\gamma)}{\beta}}
        {\fullblank}{\fullblank}{\fullblank}
    =\vm{\dpr{\alpha}{\beta}}{\fullblank}{\fullblank}{\fullblank}.
\end{displaymath}
But $\dpr{\alpha}{\beta}=0$, as $\det{x}=1$.
\end{proof}

\begin{lemma}\label{Lm:ConjugationsPermute}
Let $x\in G\le M$ be an element of order $3$, $G\cong S_3$, and $M$ a Moufang
loop. Then $T_x\in\aut{M}$ permutes the involutions of $G$.
\end{lemma}
\begin{proof}
Let $G=\span{x,\,y}\cong S_3$ with $\order{y}=2$. The remaining two
involutions of $G$ are $xy$ and $yx$. Then $T_x(y)=x^{-1}yx=xy$, $T_x(xy)=yx$,
and $T_x(yx)=y$. This can be seen from any presentation of $G$, or easily via
the natural representation of $S_3$ with $x=(1,\,2,\,3)$, $y=(1,\,2)\in S_3$.
\end{proof}

We are going to show that $\aut{\paige{2}}$ acts transitively on the subgroups
of $\paige{2}$ isomorphic to $\cyclic{2}$ --- the \emph{copies} of $\cyclic{2}$
in $\paige{2}$. Let
\begin{displaymath}
    x_0=\inv{111}{111}
\end{displaymath}
be the canonical involution. For a vector $\alpha$, let $\weight{\alpha}$ be
the \emph{weight} of $\alpha$, i.e., the number of nonzero coordinates of
$\alpha$.

\begin{prop}\label{Pr:TransitiveC2}
The group $\aut{\paige{2}}$ acts transitively on the copies of $\cyclic{2}$ in
$\paige{2}$.
\end{prop}
\begin{proof}
Let $x=\iinv{\alpha}{\beta}{a}$ be an involution. We transform $x$ into $x_0$.
By Proposition \ref{Pr:Involutions}(iii), $x$ is contained in some $G\cong
S_3$. By Lemma \ref{Lm:ConjugationsPermute} and Proposition
\ref{Pr:Involutions}(iv), we may assume that $a=0$.

Let $r=\weight{\alpha}$, $s=\weight{\beta}$. Using the automorphism $\dswitch$
from Lemma \ref{Lm:DiagonalSwitch} we can assume that $r\ge s$. We now
transform $x$ into $x'$ so that $x'=x_0$ or $x'=x_1$ or $\span{x',\,x_0}\cong
S_3$ or $\span{x',\,x_1}\cong S_3$, where $x_1=\inv{100}{100}$.

If $r\not\equiv s\pmod 2$, then $\span{x,\,x_0}\cong S_3$, by Proposition
\ref{Pr:Involutions}(ii). Suppose that $r\equiv s$. Every permutation of
coordinates can be made into an automorphism of $\paige{2}$, by Lemma
\ref{Lm:ShuffleCoords}. The involution $x_0$ is invariant under all
permutations. Since $a=0$, we must have $s>0$, and thus $(r,\,s)=(2,\,2)$,
$(1,\,1)$, $(3,\,1)$ or $(3,\,3)$. If $(r,\,s)=(2,\,2)$, transform $x$ into
$x'=\inv{110}{011}$, and note that $\span{x',\,x_1}\cong S_3$. If
$(r,\,s)=(1,\,1)$, transform $x$ into $x'=x_1$. If $(r,\,s)=(3,\,1)$,
transform $x$ into $x'=\inv{111}{001}$. Once again, $\span{x',\,x_1}\cong
S_3$. Finally, if $(r,\,s)=(3,\,3)$, we have $x=x'=x_0$.

Now, when $\span{x',\,x_i}\cong S_3$ for some $i\in\{0,\,1\}$, we can permute
the involutions of $\span{x',\,x_i}$ so that $x'$ is mapped onto $x_i$, by
Lemma \ref{Lm:ConjugationsPermute}.

It remains to show how to transform $x_1$ into $x_0$. For that matter, consider
the element $y=\tri{001}{101}{1}$, and check that $x_0=T_y(x_1)$.
\end{proof}

Let us now continue by showing that there are at most two orbits of
transitivity for copies of $V_4$. (In fact, there are exactly two orbits but
we will not need this fact.)

Put
\begin{displaymath}
    u_0=\inv{000}{110},\;\;\;
    u_1=\inv{001}{001},\;\;\;
    u_2=\inv{100}{010}.
\end{displaymath}

\begin{lemma}\label{Lm:AtMostTwo}
Let $V_4\cong\span{u,\,v}\le\paige{2}$. There is $f\in\aut{\paige{2}}$ such
that $f(u)=x_0$ and $f(v)$ is one of the two elements $u_1$, $u_2$.
\end{lemma}
\begin{proof}
By Proposition \ref{Pr:TransitiveC2}, we may assume that $u=x_0$. Write
$v=\inv{\alpha}{\beta}$, $r=\weight{\alpha}$, $s=\weight{\beta}$. We have
$r\equiv s$, else $\span{u,\,v}\cong S_3$. Thanks to the automorphism
$\dswitch$, we may assume that $r\le s$. If $(r,\,s)=(0,\,2)$, transform $v$
into $u_0$; if $(r,\,s)=(1,\,1)$, into $u_1$ or $u_2$; if $(r,\,s)=(1,\,3)$,
into $u_3=\inv{001}{111}$; if $(r,\,s)=(2,\,2)$, into $u_4=\inv{110}{110}$ or
$u_5=\inv{011}{101}$.

Let $v_1=\tri{010}{110}{0}$, $v_2=\tri{001}{101}{0}$, and define
$f_1=T_{v_2^{-1}}\circ T_{v_1}$, $f_2=T_{v_1^{-1}}\circ T_{v_2}$ (compose
mappings from the right to the left). Then $f_1$, $f_2$ are automorphisms of
$\paige{2}$, and one can check directly that $f_1(x_0)=f_2(x_0)=x_0$. Moreover,
$f_1(u_4)=u_1$, $f_1(u_3)=u_2$, $f_1(u_5)=u_3$, and $f_2(u_5)=\dswitch(u_0)$.
Thus $u_4$ can be transformed into $u_1$, and each of $u_0$, $u_3$, $u_5$ into
$u_2$.
\end{proof}

\section{Main Result}

\noindent We are now ready to demonstrate that the map $h:\octo{2}\to\octo{2}$
constructed in section \ref{Sc:RE} is an extension of $g$.

\begin{prop}\label{Pr:Additivity}
Let $C$ be a composition algebra, and let $M\subseteq C$ be the set of all
elements of norm $1$. Assume that $x$, $y\in M$ are such that $x+y\in M$. Then
$g(x+y)=g(x)+g(y)$ for every $g\in\aut{M}$.
\end{prop}
\begin{proof}
If $x=y$, we have $1=\norm{x+y}=\norm{2x}=4\norm{x}=4$. Therefore the
characteristic is $3$, and $g(x)+g(x)=-g(x)=g(-x)=g(x+x)$.

Assume that $x\ne y$. By Lemma \ref{Lm:NormOne}, $\order{xy^{-1}}=3$, and so
$\order{g(x)g(y)^{-1}}=\order{g(xy^{-1})}=3$ as well. Then
$\norm{g(x)+g(y)}=1$, again by Lemma \ref{Lm:NormOne}. Consequently, we use
Lemma \ref{Lm:MultToAdd} twice to obtain $g(x)+g(y) = -g(x)g(y)^{-1}g(x) =
g(-xy^{-1}x) = g(x+y)$.
\end{proof}

We proceed to prove by induction on the number of summands that
\begin{displaymath}
    g(\sum_{i=1}^n x_i)=\sum_{i=1}^n g(x_i)
\end{displaymath}
for every $g\in\aut{\paige{2}}$ and $x_1$, $\dots$, $x_n\in\paige{2}$ such
that $x_1+\cdots +x_n\in\paige{2}$.

\begin{lemma}\label{Lm:ItIsV4}
Suppose that $x$, $y\in\paige{2}$, $x\ne y$, are such that none of
$x+\neutral$, $y+\neutral$, $x+y$ belongs to $\paige{2}$. Then
$\span{x,\,y}\cong V_4$, and there are $a$, $b\in\paige{2}$ such that
$a+b=\neutral$, and $x+a$, $y+b\in\paige{2}$.
\end{lemma}
\begin{proof}
We have $\norm{x+\neutral}=0$, i.e., $\bilin{x}{\neutral}=0-1-1=0$. Then, by
Lemma \ref{Lm:NormZero}, $x^2=(xe^{-1})^2=-\neutral=\neutral$. Similarly,
$y^2=(xy^{-1})^2=\neutral$.

Since $\span{x,\,y}\cong V_4$, we may assume that $(x,\,y)=(x_0,\,u_1)$ or
$(x,\,y)=(x_0,\,u_2)$, where $x_0$, $u_1$, $u_2$ are as in Lemma
\ref{Lm:AtMostTwo}. When $(x,\,y)=(x_0,\,u_1)$, let $a=\tri{011}{010}{1}$,
else put $a=\tri{110}{100}{1}$. In both cases, let $b=\neutral-a$, and verify
that $x+a$, $y+b\in\paige{2}$.
\end{proof}

\begin{prop}\label{Pr:VeryAdditive}
Let $x_1$, $\dots$, $x_n\in\paige{2}$ be such that $x=\sum_{i=1}^n x_i$ belongs
to $\paige{2}$. Then
\begin{displaymath}
    g\bigl(\sum_{i=1}^n x_i\bigr)=\sum_{i=1}^n g(x_i).
\end{displaymath}
\end{prop}
\begin{proof}
The case $n=1$ is trivial, and $n=2$ is just Proposition \ref{Pr:Additivity}.
Assume that $n\ge 3$ and that the Proposition holds for all $m<n$. We can
assume that at least two summands $x_i$ are different, say $x_{n-2}\ne
x_{n-1}$. Since $g(xx_n^{-1})=g(x)g(x_n)^{-1}$, we can furthermore assume that
$x_n=\neutral$. When at least one of $x_{n-2}+\neutral$, $x_{n-1}+\neutral$,
$x_{n-2}+x_{n-1}$ belongs to $\paige{2}$, we are done by the induction
hypothesis. Otherwise, Lemma \ref{Lm:ItIsV4} applies, and there are $a$,
$b\in\paige{2}$ such that $a+x_{n-2}$, $b+x_{n-1}\in\paige{2}$, and
$a+b=\neutral$. Therefore,
\begin{eqnarray*}
    g(x)&=&g(x_1+\cdots+ x_{n-3}+(x_{n-2}+a)+(x_{n-1}+b))\\
        &=&g(x_1)+\cdots+ g(x_{n-3}) +g(x_{n-2}+a)+g(x_{n-1}+b)\\
        &=&g(x_1)+\cdots+ g(x_{n-1})+g(a)+g(b)\\
        &=&g(x_1)+\cdots+ g(x_{n-1})+g(a+b),
\end{eqnarray*}
and we are through.
\end{proof}

\begin{thm}\label{Th:Extending2}
Every automorphism of $\paige{2}$ can be uniquely extended into an
automorphism of $\octo{2}$. In particular, $\aut{\paige{2}}$ is isomorphic to
$G_2(2)$.
\end{thm}
\begin{proof}
Pick $g\in\aut{\paige{2}}$. Using the basic triple for $\octo{2}$ from Lemma
\ref{Lm:DoublingTriple} construct an automorphism $h=\psi_{\octo{2}}$ of
$\octo{2}$, as in Proposition \ref{Pr:ExtendInSteps}. Then $g$, $h$ coincide on
a basis induced by the doubling triple. Every element of $\paige{2}$ is a sum
of some of the basis elements. Hence, by Proposition \ref{Pr:VeryAdditive}, $g$
and $h$ coincide on $\paige{2}$.

This extension is unique. Thus $\aut{\paige{2}}=\aut{\octo{2}}$, and
$\aut{\octo{2}}$ is isomorphic to $G_2(2)$ by a theorem of Springer and
Veldkamp.
\end{proof}

\section{Acknowledgement}

\noindent Most of this paper is extracted from the author's Ph.\ D.\ thesis
\cite{PetrPhD}. I would like to acknowledge the support of the Department of
Mathematics at Iowa State University. I also thank the Grant Agency of Charles
University for partially supporting my visit to Oxford (grant number
269/2001/B-MAT/MFF).

Shortly before this paper was accepted for publication, G\'abor P.\ Nagy and
the author proved \cite{NagyVojtechovsky}, using completely different methods,
that $\aut{\paige{k}}$ is the semidirect product $G_2(k)\rtimes\aut{k}$, for
every perfect field $k$.

\bibliographystyle{plain}

\end{document}